\documentclass[centertags]{amsart}
\usepackage{amsmath,amstext,amsthm, amssymb,amscd}

\usepackage[mathscr]{eucal}
\usepackage{mathrsfs}

\numberwithin{equation}{section}

 \def\cU{\mathscr{U}}

\def\kg{\mathfrak{g}}

\DeclareMathOperator{\Ker}{Ker}

 \DeclareMathOperator{\ind}{index}

\newcommand{\om}{\omega}

\newtheorem{thm}{Theorem}[section]
\newtheorem{lemma}[thm]{Lemma}
\newtheorem{prop}[thm]{Proposition}

\theoremstyle{definition}
\newtheorem{rem}[thm]{Remark}
\theoremstyle{definition}
\newtheorem{defn}[thm]{Definition}
\newcommand{\be}{\begin{eqnarray}}
\newcommand{\ee}{\end{eqnarray}}

\newcommand{\comment}[1]{}

\begin{document}
\title{Geometric quantization for proper actions}

\author{Varghese Mathai}
\address{Department of Mathematics, University of Adelaide,
Adelaide 5005, Australia}
\email{mathai.varghese@adelaide.edu.au}

\author{Weiping Zhang \\ \\(With an Appendix by Ulrich Bunke)}
\address{Chern Institute of Mathematics \& LPMC, Nankai University, Tianjin 300071, PR China}
\email{weiping@nankai.edu.cn}

\begin{abstract}
 We first introduce an invariant index for $G$-equivariant elliptic
 differential
 operators on a  locally compact manifold $M$ admitting a proper
 cocompact action of a locally compact group $G$. It generalizes the
 Kawasaki index for orbifolds to the case of proper cocompact actions.
Our invariant index is used to show
 that an analog of the Guillemin-Sternberg geometric quantization
 conjecture holds if $M$ is  symplectic  with a Hamiltonian action of $G$
 that is proper and cocompact. This essentially solves a
 conjecture of Hochs and Landsman.
 \end{abstract}


\maketitle

\section{Introduction} \label{s1}

The main purpose of this paper is to generalize the
Guillemin-Sternberg geometric quantization conjecture \cite{GuSt82},
which was proved in \cite{M},  to the case of noncompact spaces and
group actions. Here we will consider the framework considered by
Hochs and Landsman in \cite{HL}.

To be more precise, let $(M,\om)$ be a locally compact symplectic
manifold. Assume that there exists a Hermitian line bundle $(L,h^L)$
over $X$ carrying a Hermitian connection $\nabla^L$ such that
$\frac{\sqrt{-1}}{2\pi}(\nabla^L)^2=\omega$. Let $J$ be an almost
complex structure on $TM$ such that $\om(\cdot, J\cdot)$
 defines a Riemannian metric $g^{TM}$ on $TM$.
  Let $D^L:\Omega^{0,*}(M,L)\rightarrow \Omega^{0,*}(M,L)$
be the canonically associated Spin$^c$-Dirac operator (cf.
\cite[Section 1]{TZ98}). Let $D^L_\pm$ be the restrictions of
$D^L$ on $\Omega^{0,{{\rm even}\over{\rm odd}}}(M,L)$
respectively.

Let $G$ be a locally compact   group with Lie algebra $\kg$.
Suppose that $G$ acts on $M$ properly. The proper $G$-action
ensures that the isotropy subgroups $G_x = \{g \in G | gx = x \}$
are compact subgroups of $G$ for all $x\in M$, all the orbits $Gx
= \{gx | g \in G\} \subseteqq M$ are closed, and moreover the
space of orbits $M/G$ is Hausdorff (cf. \cite{Palais}).

We assume  that the action of $G$ on $M$ lifts on $L$. Moreover,
we assume the $G$-action preserves the above metrics and
connections  on $TM, L$, and $J$. Then $ D_{\pm}^L $ commute with
the $G$-action.

The action of $G$ on $L$ induces naturally a moment map $\mu :
M\to \kg^*$ such that for any $V\in \kg$, $s\in\Gamma(L)$, if
$V_M$ denotes the induced Killing vector field on $M$, then the
following Kostant formula for Lie derivative holds,
\begin{align}\label{1.0}
L_{V_M}s=\nabla_{V_M}^Ls-2\pi\sqrt{-1}\langle \mu,V\rangle s.
\end{align}

 We make the assumption that $0\in \kg^*$ is a regular value
of $\mu$ and that $G$ acts freely on $\mu^{-1}(0)$. Then the
Marsden-Weinstein symplectic reduction $(M_G=\mu^{-1}(0)/G,
\om_{M_G})$ is a symplectic manifold. Moreover, $(L,\nabla ^L)$
 descends to $(L_G,\nabla ^{L_G})$ over $X_G$ so that the corresponding curvature condition
$\frac{\sqrt{-1}}{2\pi} R^{L_G}= \omega_G$ holds. The
$G$-invariant almost complex structure $J$ also descends to an
almost complex structure on $TM_G$, and $h^L$, $ g^{TM}$ descend
to $h^{L_G}$, $g^{TM_G}$ respectively. Let  $D^{L_G}$ denote the
corresponding Spin$^c$-Dirac operator on $M_G$.

Following \cite{HL}, we make the assumption that the quotient
space $M/G$  is compact, that is, the $G$-action on $M$ is
cocompact. Then $M_G=\mu^{-1}(0)/G$ is also compact.

In this paper, we will first define in Section \ref{s2} what we
call the $G$-invariant index associated to $D^L_+$, denoted by
${\rm ind}_G(D^L_+)$, which generalizes the usual definition in
the compact case to the non-compact case.

\comment{\begin{thm} \label{t1.1}    If the   Lie algebra  $\kg$
admits an ${\rm Ad}\, G$-invariant metric, then one has,
\begin{align}\label{1.1}\dim\left(\Ker D_+^L\right)^G-\dim\left(\Ker D_-^L\right)^G
={\rm ind}\left( D_+^{L_G}\right).
\end{align}
\end{thm}}

For any positive integer $p$, let $L^p$ denote the $p$-th tensor
power of $L$. Then it admits the canonically induced $G$-action as
well as the $G$-invariant Hermitian metric and connection.

We can now state the main result of this paper, which might be
thought of as a ``quantization commutes with reduction" result, in
the sense  of a conjecture of Hochs and Landsman \cite[Conjecture
1.1]{HL}, as follows.

\begin{thm}\label{t1.3}   In the general case where $G$ is merely assumed
 to be locally compact,
there exists $p_0>0$ such that for any integer $p\geq p_0$,
\begin{align}\label{1.2}{\rm ind}_G\left(D_+^{L^p}\right)
={\rm ind}\left( D_+^{L_G^p}\right).
\end{align}
Moreover, if ${\bf g}^*$ admits an ${\rm Ad}_G$-invariant metric,
then one can take $p=1$ in (\ref{1.2}).
\end{thm}

\begin{rem}\label{t1.2} In the special case where $G$ (and thus $M$) are compact, Theorem \ref{t1.3} is
the Guillemin-Sternberg geometric quantization conjecture
\cite{GuSt82} first proved by Meinrenken in \cite{M} (See  the
excellent survey of Vergne \cite{V} for further related works). In
the special case where $G$ is noncompact and admits a normal
discrete subgroup $\Gamma$ such that $\Gamma$ acts on $M$ freely
and that $G/\Gamma$ is compact,\footnote{In this case there is an
${\rm Ad}_G$ invariant metric on ${\bf g}^*$.}
   Theorem \ref{t1.3} is closely related to
    \cite[Theorem 1.2]{HL}.
 While in the special case where $G$ is semisimple and acts on
$G\times_K N$, with $K$ the maximal compact subgroup of $G$ acting
 Hamiltonianly  on a compact symplectic manifold $N$, Theorem
\ref{t1.3} should be closely related to the quantization formula of
Hochs obtained in \cite{H}. In fact, in the formulas of Hochs and
Hochs-Landsman, the left hand side of the quantization formula is
interpreted by using noncommutative $K$-theories. Thus, in combining
with our result, in the case considered by them, our invariant index
admits the noncommutative $K$-theoretic interpretation. In fact, in
the Appendix to this paper, Bunke establishes such an
interpretation. Combining  with Kasparov's index theorem
\cite{Kasp83}, one   gets a topological counterpart to our analytic
index.
\end{rem}

Since here one is dealing with noncompact group actions on
noncompact spaces, one can not apply the Atiyah-Bott-Segal-Singer
equivariant index theorem directly as in \cite{M} to prove Theorem
\ref{t1.3}. Instead, we will generalize the analytic proof of the
Guillemin-Sternberg conjecture due to Tian and Zhang \cite{TZ98} to
the current situation.

On the other hand, one can show that the finiteness of $\dim(\Ker
D_\pm^L)^G$, the dimensions of the $G$-invariant subspaces of
$\Ker(D_\pm^L)$,  holds for any equivariant Dirac type operator on
spaces with proper cocompact actions. Moreover,   we will show
that the vanishing properties of the half kernel of Spin$^c$-Dirac
operators due to Braverman \cite{B}, Borthwick-Uribe \cite{BU} and
Ma-Marinescu \cite{MM} still hold for $D^{L^p}$ when $p>0$ is
large. Combining with Theorem \ref{t1.3}, one gets

\begin{thm}\label{t1.6} There exists $p_1>0$ such that for any
integer $p\geq p_1$, one has
\begin{align}\label{1.3}
\dim\left(\Ker D^{L^p}_+\right)^G=\dim\left(\Ker
D_+^{L_G^p}\right),
\end{align}
$$\dim\left(\Ker D^{L^p}_-\right)^G=\dim\left(\Ker
D_-^{L_G^p}\right)=0.$$
\end{thm}

\begin{rem}\label{t1.4} Just as in \cite{TZ98}, one can twist
$L^p$ by an arbitrary $G$-equivariant vector bundle $F$ over $M$
carrying a $G$-invariant Hermitian metric and a $G$-invariant
Hermitian connection (that is, replace $L^p$ by $L^p\otimes F$).
Then all the results of this paper still hold when $p>0$ is large
enough. In particular, if we chose $F=\Lambda^{*,0}(T^*M)$, then
we get a non-compact quantization formula for the signature
quantization considered in \cite{TZ1} and \cite{GSW}.
\end{rem}

\begin{rem}\label{t1.5} Our method also works for the case where
the action of $G$ on $\mu^{-1}(0)$ is not free. Then
$M_G=\mu^{-1}(0)/G$ is an orbifold, and (\ref{1.2}) and
(\ref{1.3}) still hold if we replace the righthand sides by the
corresponding orbifold indices. We leave this to the interested
reader.
\end{rem}

The rest of this paper is organized as follows. In \S\ref{s2},
for any locally compact manifold $M$ admitting a proper
cocompact action of a locally compact group $G$, and any
$G$-equivariant Dirac type operator $D$ on $M$, we introduce what
we call the $G$-invariant index ${\rm ind}_G(D)$. In
\S\ref{s3}, we generalize the analytic techniques developed in
\cite{TZ98} to give a proof of Theorem \ref{t1.3}. In
\S\ref{s4}, we prove the vanishing properties of the $G$-invariant
part of the half kernel of the equivariant Dirac operator and as a
consequence, get Theorem \ref{t1.6}.  Finally, in \S5, we
consider some examples and applications of our main results.

\section{The invariant index for proper cocompact actions} \label{s2}

Let $M$ be a locally compact manifold. Let $G$ be a locally
compact group. Let $dg$ be the left invariant Haar measure on $G$.

We make the assumption that $G$ acts on $M$ properly and
cocompactly,  where by proper action we mean that the following
map
\begin{align}\label{2.1}
G\times M\rightarrow  M\times M,\ \ \ \ \ \ (g,x)\mapsto (x,gx)
\end{align} is proper (that is, the inverse image of a compact
subset is compact), while by cocompact we mean that the quotient
$M/G$ is a compact space.

One of the basic properties for such an action is that there
exists a smooth, non-negative, compactly supported cut-off
function $c$ on $M$ such that
\begin{align}\label{2.2}\int_Gc(g^{-1}x)^2dg=1
\end{align}
 for any $x\in M$ (cf.
\cite[Section 7.2.4, Proposition 8]{Bo}). This cut-off function
allows one to get $G$-invariant objects from the originally not
necessarily $G$-invariant ones.

As an example, for any Riemannian metric $g^{TM}$ on $TM$ one gets
an ``averaged" $G$-invariant metric
\begin{align}\label{2.3}
\int_Gc(h^{-1}x)^2(h^*g^{TM})(x)dh \end{align}
 on $TM$.

 From now on, we assume that $g^{TM}$ is $G$-invariant.

 Let ${\rm cl}(TM)$ be the Clifford algebra bundle associated to
 $(TM,g^{TM})$. Then it admits a naturally  induced $G$-action as well as a $G$-invariant Hermitian
 metric $g^{{\rm cl}(TM)}$.

 Let $E$ be a complex vector bundle over $M$ such that $E=E_+\oplus E_-$ is a ${\bf Z}_2$-graded
 ${\rm cl}(TM)$-module and that
 it admits a
 $G$-action, preserving the ${\bf Z}_2$-grading of $E$,
  lifted from the action of $G$ on $M$. Let $g^E$ be a ${\bf Z}_2$-graded
 $G$-invariant Hermitian metric on $E$, let $\nabla^E$ be a ${\bf Z}_2$-graded
 $G$-invariant Hermitian connection on $E$. The averaging procedure
 similar to that in (\ref{2.3}) guarantees the existence of $g^E$
 and $\nabla^E$.

Let $e_1,\cdots,e_{\dim M}$ be an oriented orthonormal basis of
$TM$.  We define a Dirac type operator $D^E$ to be the operator
acting on $\Gamma(E)$,
\begin{align}\label{2.4}
D^E=\sum_{i=1}^{\dim M}c(e_i)\nabla^E_{e_i}+A:\Gamma(E)\rightarrow
\Gamma(E),
\end{align}
 where $A\in \Gamma({\rm End}(E))$   exchanges $E_\pm$.

 We make the assumption that $D^E$ is $G$-equivariant.

 Let $\Gamma(E)$  carry the
 natural inner product such that
 for any $s,\, s'\in \Gamma(E)$ with compact supports,
\begin{align}\label{2.5}
\left\langle s,s'\right\rangle =\int_M\left\langle
s(x),s'(x)\right\rangle_E dx.
\end{align}
Let $\|\cdot\|_0$ denote the associated $L^2$-norm. Let $L^2(M,E)$
denote the  completion of $\Gamma(E)$ with respect to the inner
product $\|\cdot\|_0$.

Since $M/G$ is compact, there exists a compact subset $Y$ of $M$
such that $G(Y)=M$ (cf. \cite[Lemma 2.3]{P}).

Let $U$, $U'$ be two open subsets of $M$ such that $Y\subset U$
and that the closures $\overline{U}$ and $\overline{U'}$ are both
compact in $M$, and that $\overline{U}\subset U'$. The existence
of $U$, $U'$ is clear.

Then it is easy to construct a smooth function $f:M\rightarrow
[0,1]$ such that $f|_U=1$ and ${\rm Supp}(f)\subset
U'$. \footnote{With this function, one can construct the cut-off
function $c$ by $c(x)={f(x)\over(\int_Gf(g^{-1}x)^2dg)^{1/2}}$ for
any $x\in M$.}

We now consider the space $\Gamma(E)^G$, the subspace of
$G$-invariant sections of $\Gamma(E)$.

By using the property that $G(Y)=M$, it is easy to see that there
exists a positive constant $C>0$ such that for any
$s\in\Gamma(E)^G$,
\begin{align}\label{2.6}
\|s\|_{U,0}\leq \|fs\|_0\leq \|s\|_{U',0}\leq C\|s\|_{U,0},
\end{align} where for $V=U$ or $U'$,
\begin{align}\label{2.7}
\|s\|_{V,0}^2=\int_V\left\langle s(x),s(x)\right\rangle_E dx.
\end{align}

Let ${\bf H}^0_f(M,E)^G$ be the completion of the space $\{fs:
s\in \Gamma(E)^G\}$ under the norm $\|\cdot\|_0$ associated to the
inner product (\ref{2.5}). Let ${\bf H}^1_f(M,E)^G$ be the
completion of $\{fs: s\in \Gamma(E)^G\}$ under a (fixed)
$G$-invariant  first Sobolev norm associated to the inner product
(\ref{2.5}).

Let $P_f$ be the orthogonal projection from $L^2(M,E)$ to its
subspace ${\bf H}^0_f(M,E)^G$.

It is clear that $P_fD^E$ maps an element of ${\bf H}^1_f(M,E)^G$
into ${\bf H}^0_f(M,E)^G$.

\begin{prop}\label{t2.1} The induced operator
\begin{align}\label{2.8}
P_fD^E:{\bf H}^1_f(M,E)^G\rightarrow {\bf H}^0_f(M,E)^G
\end{align}
is a Fredholm operator.

\end{prop}

{\it Proof}. For any $s\in \Gamma(E)^G$, by (\ref{2.4}), one has
\begin{align}\label{2.9}
D^E(fs)=fD^Es+c(df)s,
\end{align}
where we identify the one form $df$ with its metric dual $(df)^*$.

Since $D^E$ is $G$-equivariant, it is clear that $P_f(fD^Es)=fD^Es$,
while in view of (\ref{2.6}),
\begin{align}\label{2.10}
\left\|P_f(c(df)s)\right\|_0\leq \|c(df)s\|_0\leq C_1\|s\|_{U,0}
\end{align}
for some constant $C_1>0$.

Thus, one has, by also proceeding as in (\ref{2.6}),
\begin{align}\label{2.11}
\left\|P_fD^E(f s)\right\|_0\geq \|fD^Es\|_0- C_1\|s\|_{U,0} \geq
C_2\|s\|_{U,1}-C_3\|s\|_{U,0}
\end{align}
for some constants $C_2,\,C_3>0$.

On the other hand, by (\ref{2.9}) and by proceeding as in
(\ref{2.6}), one verifies that
\begin{align}\label{2.12}
 \|fs\|_1\leq
 C_4\left(\|fs\|_0+\left\|D^E(fs)\right\|_0\right)\leq
 C_5\|s\|_{U,0}+C_6\|s\|_{U,1}.
\end{align} for some constants $C_4,\,C_5,\,C_6>0$.

From (\ref{2.6}), (\ref{2.11}) and (\ref{2.12}), one gets
\begin{align}\label{2.13}
\left\|P_fD^E(f s)\right\|_0\geq C_7\| fs\|_{1}-C_8\| fs\|_{0},
\end{align}
for some constants $C_7,\,C_8>0$.

Since $f$ is of compact support, from the G\"arding type
inequality (\ref{2.13}) one sees that $P_fD^E$ is a Fredholm
operator.
 \ \ Q.E.D.

\begin{rem}\label{t2.2} Besides the Fredholm property in
Proposition \ref{t2.1},   the following self-adjoint property also
holds:  for any $s,\,s'\in \Gamma(E)^G$, one has
\begin{align}\label{2.14}
\left\langle P_fD^E(fs),fs'\right\rangle=\left\langle
fs,P_fD^E\left(fs'\right)\right\rangle,
\end{align}if  $D^E$ is   formally self-adjoint.
 \end{rem}

\begin{rem}\label{t2.3}  If
$(\widetilde{U},\,\widetilde{U}',\,\widetilde{f})$ is another
triple of open subsets  and the cut-off function as above, then by
taking the deformation $f_t=(1-t)f+t\widetilde{f}$, one gets
easily a continuous family of Fredholm operators $P_{f_t}D^E$.
\end{rem}

Now let $D_\pm^E:\Gamma(E_\pm)\rightarrow \Gamma(E_\mp)$ be the
restrictions of $D^E$ on $\Gamma(E_\pm)$ respectively.

Then by Proposition \ref{t2.1} and (\ref{2.14}), the induced
operator $P_fD^E_+:{\bf H}^1_f(M,E_+)^G\rightarrow {\bf
H}^0_f(M,E_-)^G$ is Fredholm. Moreover its index, ${\rm
ind}(P_fD^E_+)$, does not depend of the choice of $f$, in view of
Remark \ref{t2.3}. Similarly, it is also easy to see that this
index does not depend on the choices of $G$-invariant metrics and
connections involved.

\begin{defn}\label{t2.5}
We call ${\rm ind}(P_fD^E_+)$ defined above the $G$-invariant
index associated to $D^E_+$ and denote it by ${\rm ind}_G(D^E_+)$.
\end{defn}

\begin{rem}\label{t2.6} If $M$ is compact and $D^E$ is formally self-adjoint, then one can take
$f\equiv 1$, so that one has
\begin{align}\label{2.17}
{\rm ind}_G\left(D^E_+\right)=\dim\left(\left(\ker D^E_+
\right)^G\right)-\dim\left(\left(\ker D^E_- \right)^G\right) .
\end{align}
\end{rem}

\begin{rem}\label{t2.4}   For any $s\in \Gamma(E)^G$, it is clear
that $fs$ and $s$ are determined by each other. That is, if $s,\,
s'\in \Gamma(E)^G$ are such that $fs=fs'$, then $s=s'$ as $f\equiv
1$ on $Y$. Moreover, by (\ref{2.9}), (\ref{2.10}) and (\ref{2.13}),
one sees easily that the induced operator $\widetilde{D}^E_f:{\bf
H}^1_f(M,E)^G\rightarrow {\bf H}^0_f(M,E)^G$ such that
\begin{align}\label{2.15}
\widetilde{D}^E_f (fs):=fD^E(s)
\end{align}
is a Fredholm operator. Thus one gets that
\begin{align}\label{2.16}
\dim\left(\ker\left(
\left.D^E\right|_{\Gamma(E)^G}\right)\right)=\dim\left(\left(\ker
D^E \right)^G\right)<+\infty.
\end{align}
\end{rem}

In fact, when $G$ is unimodular, we can further identify  ${\rm
ind}_G(D^E_+) $ as follows.

\begin{thm}\label{l2.7} If $G$ is unimodular, then with the notation above, one has
\begin{align}\label{a.1}
{\rm ind}\left(P_fD^E_+\right)=\dim\left(\left(\ker D^E_+
\right)^G\right)-\dim\left(\left(\ker D^E_- \right)^G\right) ,
\end{align}
whenever $D^E$ is formally self-adjoint.
\end{thm}

{\it Proof}. We use the cut-off function $f$, and set for any $x\in M$,
\begin{align}\label{a.2}
c(x)={f(x)\over\displaystyle \left(
\int_Gf(g^{-1}x)^2dg\right)^{1/2}}.
\end{align} Then
\begin{align}\label{a.3}
\int_G c(g^{-1}x)^2dg=1
\end{align} for any $x\in M$.
Let $\alpha$ denote the positive  $G$-invariant function on $M$
defined by
\begin{align}\label{a.4}
\alpha(x)=\left(\int_Gf(g^{-1}x)^2dg\right)^{1/2}.
\end{align}

Let ${\bf H}^0_c(M,E)^G$ be the $L^2$ completion of the space
$\{c s:s\in \Gamma(E)^G\}$. Let ${\bf H}^1_c(M,E)^G$ be the
corresponding first Sobolev space associated to a (fixed)
$G$-invariant first Sobolev norm.

Let $\beta_f:{\bf H}^0_f(M,E)^G\rightarrow {\bf H}^0_c(M,E)^G$
be the isomorphism such that for any $s\in \Gamma(E)^G$,
\begin{align}\label{a.5}
\beta_f: fs\longmapsto f {s\over \alpha}=c s.
\end{align}

Let $P_c$ be the orthogonal projection from $L^2(M,E)$ onto
${\bf H}^0_c(M,E)^G$.

For any $s\in \Gamma(E)^G$, one verifies that
\begin{align}\label{a.6}
\beta_f\left(P_f D^E(fs)\right)= c D^Es+{1\over
\alpha}P_f(c(df)s)= P_c D^E(c s)-P_c(c(dc)s)+{1\over
\alpha}P_f(c(df)s).
\end{align}

Since $P_fD^E$ is Fredholm, from (\ref{a.6}), one sees easily that
$P_c D^E$ is a Fredholm operator and
\begin{align}\label{a.7}
{\rm ind}\left(P_f D^E_+ \right)= {\rm ind}\left( P_c
D^E_+\right) .
\end{align}

Indeed, what Bunke does in his Appendix is to give a
$KK$-theoretic interpretation of the right hand side of
(\ref{a.7}).

Moreover,  Bunke actually writes out explicitly the projection
$P_c$ when $G$ is unimodular.
 According to Bunke, when $G$ is unimodular, for any $\mu\in \Gamma (E)$ with compact
 support, one has (cf.   Appendix D, where $P_c$ here is exactly $Q_\epsilon$ there with $\epsilon =1$)
\begin{align}\label{a.8}
\left(P_c \mu\right)(x) =c (x) \int_G
c(g^{-1}x)\left(g^*\mu\right)(x) dg.
\end{align}

From (\ref{a.3}) and (\ref{a.8}), one computes that for any
$s\in\Gamma(E)^G$, one has at $x\in M$ that
\begin{align}\label{a.9}
 P_c (c(dc)s)   =c   \int_G
c(g^{-1}x)\left(g^*(c(dc)s)\right)(x) dg={1\over 2}c \cdot
c\left(d\int_G c(g^{-1}x)^2dg\right)s(x)=0.
\end{align}

From (\ref{a.9}), one gets that for any $s\in \Gamma(E)^G$,
\begin{align}\label{a.10}
 P_c D^E(c s)=c D^Es+P_c (c(dc)s)= c D^Es  .
\end{align}

When $D^E$ is formally self-adjoint, from (\ref{a.10}) one gets
immediately that
\begin{align}\label{a.11}
 \dim \left(\Ker P_c D^E_\pm\right)=\dim\left(\left(\Ker   D^E_\pm\right)^G\right)  .
\end{align}

Combining with (\ref{a.7}), one gets
\begin{align}\label{a.12}
{\rm ind}_G\left(D^E\right)={\rm ind}\left(P_c D^E\right)=
\dim\left(\left(\ker D^E_+
\right)^G\right)-\dim\left(\left(\ker D^E_- \right)^G\right) ,
\end{align}
which is exactly (\ref{a.1}).  \ Q.E.D.

\medskip
 \begin{rem}\label{t2.9}
 Remark \ref{t2.6} and Theorem \ref{l2.7} fully justify   the
 term ``$G$-invariant index'' in Definition \ref{t2.5}. Moreover,
 by (\ref{a.11}), one sees that $\ker (P_cD^E)$ consists of smooth
 elements.
 \end{rem}

\begin{rem}\label{t2.88}
 When $G$ is non-unimodular, Theorem  \ref{l2.7} still holds if one
 inserts the modular factor $\delta$ (with $dg^{-1}=\delta(g) dg$) in the right hand side of (\ref{a.1}) as in Bunke's Appendix. We
 leave it to the interested reader.
 \end{rem}

We now consider the very special case where $G$ acts on $M$
freely, but we no longer assume that $G$ is unimordular (thus we
no longer have Theorem \ref{l2.7}). Then $M/G$ is a compact
manifold, while $E$ descends to a Hermitian vector bundle $E_G$
over $M/G$ carrying an induced Hermitian connection.

Let $D_+^{E_G}:\Gamma({E_{+,G}})\rightarrow \Gamma({E_{-,G}})$ be
associated Dirac operator.

\begin{prop}\label{t2.7} The following identity holds,
\begin{align}\label{2.18}
{\rm ind}_G\left(D^E_+\right)={\rm ind}\left(D_+^{E_G}\right).
\end{align}
\end{prop}

{\it Proof.} For any $s\in\Gamma(E_G)$, let $\widetilde{s}\in
\Gamma(E)^G$ be its canonical lift.

It is clear that the map $\Gamma(E_G)\rightarrow  \{
f\Gamma(E)^G\}$ such that $s\rightarrow f\widetilde{s}$ extends
canonically to a bounded linear isomorphism
$\alpha_f:L^2(E_G)\rightarrow {\bf H}^0_f(M,E)^G$.

\begin{lemma}\label{t2.8}There exists a constant $C>0$ such that
for any $s\in \Gamma(E_G)$, one has
\begin{align}\label{2.19}
 \left\|\left(\alpha_f\right)^{-1}P_fD^E\alpha_f(s)-D^{E_G}s\right\|\leq
 C\| s\|.
\end{align}
\end{lemma}

{\it Proof}. From (\ref{2.9}), one has
\begin{align}\label{2.20}
 D^E\alpha_fs=D^E({f}\widetilde{s})=fD^E\widetilde{s}+c(df)\widetilde{s}.
\end{align}

Now, on the principal fibre bundle $G\rightarrow M\rightarrow M/G$,
the vertical  directions are generated by elements in the Lie
algebra ${\bf g}$. Thus, for any $X\in T^VM$ such that $|X|=1$,
since $\widetilde{s}$ is $G$-invariant so that
$L^E_X\widetilde{s}=0$, where $L^E_X $ is the Lie derivative along
$X$ on $E$, one has
\begin{align}\label{2.21}
\left\| f\left(\nabla^E_X\widetilde{s}\right)\right\|=\left\|
f\left(\nabla^E_X\widetilde{s}-L^E_X\widetilde{s}\right)\right\|\leq
C_1\| f\widetilde{s}\|
\end{align}
for a (fixed) positive constant $C_1>0$.

From    (\ref{2.21}), one verifies easily that
\begin{align}\label{2.22}
 \left\|\left(\alpha_f\right)^{-1}\left(fD^E\widetilde{s}\right)-D^{E_G}s\right\|\leq
 C_2\| s\|
\end{align}
for a (fixed) positive constant $C_2>0$.

From (\ref{2.22}), one sees immediately that
\begin{align}\label{2.23}
\left\|\left(\alpha_f\right)^{-1}\left(P_f\left(fD^E\widetilde{s}\right)\right)-D^{E_G}s\right\|\leq
 C_2\| s\|.
 \end{align}

 On the other hand, one finds easily that there are positive
 constants $C_3>0$, $C_4>0$, $C_5>0$ such that
\begin{align}\label{2.24}
\left\|\left(\alpha_f\right)^{-1}P_f\left(c(df)\widetilde{s}\right)\right\|\leq
 C_3\left\| c(df)\widetilde{s}\right\|\leq
 C_4\left\|f\widetilde{s}\right\|\leq C_5\| s\|.
 \end{align}

 From (\ref{2.20}), (\ref{2.23}) and (\ref{2.24}), one gets
 (\ref{2.19}).\  \ Q.E.D.

 $\ $

 We now return to the proof of (\ref{2.18}).

 By restricting (\ref{2.19}) to $\Gamma(E_{+,G})$, one deduces
  that
  \begin{align}\label{2.25}
{\rm ind}\left(D^{E_G}_+\right)={\rm
ind}\left(\left(\alpha_f\right)^{-1}P_fD^E_+\alpha_f\right)={\rm
ind}\left( P_fD^E_+ \right),
 \end{align}
 which is exactly (\ref{2.18}).\ \ Q.E.D.

\begin{rem}\label{t2.10} It is easy to see that for
  Proposition \ref{t2.7} to hold, one  need only to assume that (the
  $G$-equivariant operator)
  $D^E_+$ is $G$-transversally elliptic.
 \end{rem}

\section{The geometric quantization formula for proper
actions}\label{s3}
 We now turn back  to the situation as in \S\ref{s1}. In this
case, $E=\Lambda^{0,*}(T^*M)\otimes L^p$ and
$\Omega^{0,*}(M,L^p)=\Gamma(\Lambda^{0,*}(T^*M)\otimes L^p).$

Let $D^{L^p}$ denote the corresponding Spin$^c$-Dirac operator (cf.
\cite[Section 1]{TZ98}). Clearly, $D^{L^p}$ is $G$-equivariant.

Let $\Omega^{0,*}(M,L^p)^G$ be the subspace of $G$-invariant
sections of $\Omega^{0,*}(M,L^p)$. It is clear that any section in
$\Omega^{0,*}(M,L^p)^G$ is determined by its restriction to $Y$.

Note that since in the general case where $G$ is assumed to be
only locally compact, there might not be any ${\rm
Ad}_G$-invariant metric on $\kg$.

Choose any    metric on $\kg^*$.  Let $h_1,\, \cdots,\, h_{\dim
G}$ be an orthonormal basis of $\kg^*$. Denote by $V_i$ the
Killing vector field  on $X$ generated by the dual of $h_i$
($1\leq i\leq \dim G$). The point here is that the function
\begin{align}\label{3.1}{\mathcal
H}=\|\mu\|^2=\sum_{i=1}^{\dim G}\mu_i^2
\end{align}
might not be $G$-invariant, thus  the associated Hamiltonian
vector field $X^{\mathcal H}$ might not be $G$-invariant. We first
construct an invariant one out of it.

Recall that the cut-off function $c$ has been defined in
(\ref{2.2}).

Let
\begin{align}\label{3.2}
X_G^{\mathcal H}=\int_Gc(g^{-1}x)^2 X_g^{\mathcal H}dg
\end{align}
denote the averaged {\it $G$-invariant} vector field on $M$, where
$X_g^{\mathcal H}$ denotes  the pullback of $X^{\mathcal H}$ by
$g \in G$.

For any $T\geq 0$, set
\begin{align}\label{3.3}D^{L^p}_{T}=D^{L^p}+{\sqrt{-1}T\over 2}c\left(X_G^{\mathcal
H}\right).\end{align}
 Then it is $G$-equivariant. Moreover, it is a formally self-adjoint Dirac type operator
 in the sense of (\ref{2.4}) and thus the results in \S 2
 apply here to $D^{L^p}_T$.

From   (\ref{3.2}) and the fact that
\begin{align}\label{3.4}X^{\mathcal H}=2\sum_{i=1}^{\dim
G}\mu_iV_i\end{align} (cf. \cite[(1.19)]{TZ98}), it is clear that
at any $x\in M$, $X_G^{\mathcal H}$ lies in $T_x(Gx)$. Moreover,
by (\ref{1.0}) and (\ref{3.4}),
 one verifies that
for any $s\in \Omega^{0,*}(M,L^p)^G$,
\begin{align}\label{3.5}
 \nabla^{\Lambda^{0,*}(T^*M)\otimes L^p}_{X_G^{\mathcal
H}} s=\left(A\otimes {\rm Id}_{L^p} +4p\pi\sqrt{-1}{\mathcal
H}_G(x)\right)s,
\end{align}
where \begin{align}\label{3.6} {\mathcal
H}_G=\int_Gc(g^{-1}x)^2{\mathcal H}(g^{-1}x)dg\end{align} and
\begin{align}\label{3.7}
A= \nabla^{\Lambda^{0,*}(T^*M)}_{X_G^{\mathcal H}}
-\int_Gc(g^{-1}x)^2\left(2\sum_{i=1}^{\dim
G}\mu_i(g^{-1}x)L^{\Lambda^{0,*}(T^*M)}_{g^*V_i}\right)dg
\end{align}
 is of order zero
and does not involve $p$.

Let   $U'$ be constructed as in \S 2. Let   $W$ be any open
neighborhood of $\mu^{-1}(0)$ in $M$. We first show that the
following analogue of \cite[Theorem 2.1]{TZ98} holds.

\begin{prop}\label{t3.1} There exists  $p_0\geq 1$ such that for any integer
 $p\geq p_0$, there exist  $C>0, \, b>0$
verifying the following property:  for any  $T\geq 1$ and $s\in
\Omega^{0,*}(M,L^p)^G$ with ${\rm Supp}(s)\cap
\overline{U'}\subset \overline{U'}\setminus W $, one has
\begin{align}\label{3.8}
\left\|P_fD^{L^p}_T(fs)\right\|_0^2\geq
C\left(\|fs\|^2_1+(T-b)\|fs\|_0^2\right).
\end{align}
Moreover, if ${\bf g}^*$ admits an ${\rm Ad}_G$-invariant metric,
then one can take $p_0=1$.
\end{prop}

{\it Proof}. One computes first that
\begin{align}\label{3.9}
\left\|P_fD^{L^p}_T(fs)\right\|_0=\left\|P_f\left(fD^{L^p}_Ts+c(df)s\right)\right\|_0\geq
\left\| fD^{L^p}_Ts\right\|_0-C_1\|fs\|_0
\end{align}
for some constant $C_1>0$.

On the other hand, one has
\begin{align}\label{3.10}
\left\| fD^{L^p}_Ts\right\|_0=\left\|
D^{L^p}_T(fs)-c(df)s\right\|_0\geq \left\|
D^{L^p}_T(fs)\right\|_0-C_2\|fs\|_0
\end{align}
for some constant $C_2>0$.

 From (\ref{3.9}) and (\ref{3.10}), one sees that there exists
 $C_3>0$ such that
\begin{align}\label{3.11}
\left\|P_fD^{L^p}_T(fs)\right\|_0^2\geq {1\over 2}
\left\|D^{L^p}_T(fs)\right\|^2_0-C_3\|fs\|_0^2 .
\end{align}

Since $fs$ has compact support, one has
\begin{align}\label{3.12}
 \left\|D^{L^p}_T(fs)\right\|^2_0 =\left\langle D^{L^p}_T(fs),D^{L^p}_T(fs)\right\rangle
 =\left\langle D^{L^p,2}_T(fs), fs\right\rangle.
\end{align}

Now since $s$ is $G$-invariant, from \cite[(1.26)]{TZ98},
(\ref{3.3}), (\ref{3.5})
 and (\ref{3.7}), one computes that
\begin{multline}\label{3.13}
 D^{L^p,2}_T(fs)  =D^{L^p,2}(fs)+{\sqrt{-1}T\over 2}\sum_{j=1}^{\dim M}c(e_j)
 c\left(\nabla^{TM}_{e_j}X^{\mathcal H}_G\right)(fs)\\
 -\sqrt{-1}T\left(A\otimes {\rm Id}_{L^p}\right)(fs) +4p\pi T{\mathcal H}_Gfs
 -\sqrt{-1}TX^{\mathcal H}_G(f)s+{T^2\over 4}\left|X^{\mathcal H}_G\right|^2(fs).
\end{multline}

\begin{lemma}\label{t3.2} One has ${\mathcal
H}_G^{-1}(0)=\mu^{-1}(0).$
\end{lemma}

{\it Proof}. By the definition (\ref{3.6}) of ${\mathcal H}_G$, it
is clear that for any $x\in\mu^{-1}(0)$, ${\mathcal H}_G(x)=0$.

Conversely, for any $x\notin \mu^{-1}(0)$, by (\ref{2.2}), there
exists $g\in G$ such that    $c(g^{-1}x)\neq 0$. Thus,
$c(g^{-1}x)\mu(g^{-1}x)=c(g^{-1}x)g^*\mu(x)\neq 0$, from which and
from (\ref{3.6}) one gets that ${\mathcal H}_G(x)\neq 0$.\ \ Q.E.D.

$\ $

By Lemma \ref{t3.2}, there exists a constant $\alpha>0$ such that
\begin{align}\label{3.14}
 {\mathcal H}_G(x)\geq \alpha
\end{align}
for any $x\in \overline{U'}\setminus W $.

Clearly,
\begin{align}\label{3.15}
{\rm Re}\left(\left\langle \sqrt{-1}X^{\mathcal H}_G(f)s,
fs\right\rangle\right)=0.
\end{align}

On the other hand, since $f$ has compact support in $U'$ and ${\rm
Supp}(s)\cap \overline{U'}\subset \overline{U'}\setminus W $, it
is easy to see that there exists a constant $C_4>0$ such that
\begin{multline}\label{3.16}
{\rm Re}\left(  \left\langle{\sqrt{-1}\over 2}\sum_{j=1}^{\dim
M}c(e_j)
 c\left(\nabla^{TM}_{e_j}X^{\mathcal H}_G\right)(fs)
 -\sqrt{-1}\left(A\otimes {\rm Id}_{L^p}\right)(fs),fs\right\rangle\right)
\\  +4p\pi \left\langle{\mathcal
 H}_Gfs,fs\right\rangle \geq
 \left(4p\pi\alpha-C_4\right)\|fs\|_0^2.
\end{multline}

 From  (\ref{3.12})-(\ref{3.16}), one sees that
there exist constants $C_5,\, C_6>0$ such that
\begin{align}\label{3.17}\left\|D^{L^p}_T(fs)\right\|^2_0
 =\left\langle D^{L^p,2}_T(fs), fs\right\rangle\geq C_5\|fs\|_1^2
 -C_6\|fs\|_0^2+\left(4p\pi\alpha-C_4\right)T\|fs\|_0^2.
\end{align}

Formula (\ref{3.8}) follows from (\ref{3.11}) and (\ref{3.17}) by
taking $p_0={C_4\over \pi\alpha}$ in (\ref{3.17}).

For the remaining situation  where ${\bf g}^*$ admits an ${\rm
Ad}_G$-invariant metric, in this case, both ${\mathcal H}$ and
$X^{\mathcal H}$ are $ G$-invariant, so we are in an exactly
similar situation as in \cite[Section 2]{TZ98}.

Set
\begin{align}\label{3.18}F_T^L=D_T^{L,2}+2\sqrt{-1}T\sum_{i=1}^{\dim
G}\mu_iL_{V_i}
\end{align}
as in \cite[(1.30)]{TZ98}.

One verifies directly in this case, in view of (\ref{3.4}), that
\begin{align}\label{3.19}
 \left\langle
D_T^{L,2}(fs),fs\right\rangle=\left\langle
F_T^L(fs),fs\right\rangle -\sqrt{-1}T\left\langle X^{\mathcal
H}(f)s,fs\right\rangle.
\end{align}

Now since $f$ has compact support in $U'$ and ${\rm Supp}(s)\cap
\overline{U'}\subset \overline{U'}\setminus W $, by proceeding in
exactly the same way as in \cite[Section 2]{TZ98}, one sees that
there exist constants $C_7,\, C_8>0$ such that for any $T\geq 1$,
\begin{align}\label{3.20}
 {\rm Re}\left(\left\langle F_T^L(fs),fs\right\rangle\right)\geq
C_7\left(\|fs\|_1^2+\left(T-C_8\right)\|fs\|_0^2\right).
\end{align}

From (\ref{3.11}), (\ref{3.12}), (\ref{3.15}), (\ref{3.19}) and
(\ref{3.20}), one sees that Proposition \ref{t3.1} holds   for
$p_0=1$ in the case where ${\bf g}^*$ admits an ${\rm
Ad}_G$-invariant metric.

The proof of Proposition \ref{t3.1} is completed. \ \ Q.E.D.

$\ $

\begin{rem}\label{t3.3} Note that  in the general case where ${\bf g}^*$ does not
admit an ${\rm Ad}_G$-invariant metric,
 $X^{\mathcal H}_G$ might not be the Hamiltonian vector field
associated to $ {\mathcal H}_G$. This makes it difficult here to
get the pointwise estimates  like in \cite[Proposition 2.2]{TZ98}
at
 zeroes  of $X^{\mathcal H}_G$, and partially explains why we
need to pass to the uniform estimate for $p>0$ large.
\end{rem}

\begin{rem}\label{t3.4}  Proposition \ref{t3.1} allows us to
localize the proof of Theorem \ref{t1.3} to an arbitrarily small open
neighborhood of $\mu^{-1}(0)\cap \overline{U'}$ in
$\overline{U'}$, just as in \cite{TZ98} which relies on techniques
developed in \cite{BL91}.
\end{rem}

Indeed, for any $r>0$, let $W_r$ denote the $G$-invariant open
neighborhood of $\mu^{-1}(0)$ in $M$ such that $W_r= \{x\in M:
{\mathcal H}_G(x)<r \}. $

Since $0\in{\bf g}^*$ is a regular value of $\mu$ and $G$ acts on
$\mu^{-1}(0)$ freely, one sees easily that  when $r>0$ is small
enough, $G$ also acts on $W_r$ freely.

\begin{lemma}\label{t3.5} One has that $0\in{\bf R}$ is a non-degenerate critical
value of ${\mathcal H}_G:M\rightarrow {\bf R}$.
\end{lemma}

{\it Proof}. Let $N$ be the normal bundle to $\mu^{-1}(0)$ in $M$.
Let $g^N$ be the $G$-invariant metric on $N$ induced by the
$G$-invariant orthogonal decomposition
\begin{align}\label{3.21}
TM|_{\mu^{-1}(0)}=T\mu^{-1}(0)\oplus N,\ \ \ \ \
g^{TM|_{\mu^{-1}(0)}}=g^{T\mu^{-1}(0)}\oplus g^N.
\end{align}

Let $P^{T\mu^{-1}(0)}$, $P^N$ denote the orthogonal projections
from $TM|_{\mu^{-1}(0)}$ to $T\mu^{-1}(0)$ and $ N$ respectively
with respect to (\ref{3.21}). Let $\nabla^N$ be the connection on
$N$ defined by $\nabla^N=P^N(\nabla^{TM}|_{\mu^{-1}(0)})$ where
$\nabla^{TM}$ is the ($G$-invariant) Levi-Civita connection
associated to $g^{TM}$.

For any $x\in \mu^{-1}(0),\, Z\in N_x$, we identify $Z$ with
$\exp^{T_xM}(Z)\in M$. Since $M/G$ is compact, one verifies easily
that when $\varepsilon>0$ is small enough, the above map induces
an identification from $ N_\varepsilon=\{Z\in N:|Z|<\varepsilon\}$
to its image in $M$.

For any $y\in \mu^{-1}(0)$ with $c(y)\neq 0$, let $U_y$ be a small
enough open neighborhood of $y$ in $\mu^{-1}(0)$ such that
$c(y')\geq {1\over 2}c(y)>0$ for any $y'\in U_y$, and that there
exists $\varepsilon_y>0$ and $C_y>0$ such that ${\mathcal
H}(y',Z)\geq C_y|Z|^2$ for any $Z\in N_{y'}$ with $|y'|\leq
\varepsilon_y$. Moreover, there is an open neighborhood $G_y$ of
${e}$ in $G$ such that $c(g^{-1}y')\geq {1\over 4}c(y)$ and
${\mathcal H}(g^{-1}y',Z)\geq {1\over 2}C_y|Z|^2$ for any $g\in
G_y$ and $y'\in U_y$.
 The existence of $U_y$ is clear.

 For any $x\in \mu^{-1}(0)$, let $h\in G$ be such that $c(h^{-1}x)\neq 0$. Let
 $U_{h^{-1}x}$ be the open neighborhood of $h^{-1}x$ constructed above. Then
 $hU_{h^{-1}x}$ is an open neighborhood of $x$ such that for any
 $x'\in hU_{h^{-1}x}$ and $Z'\in N_{x'}$ with
 $|Z'|\leq\varepsilon_{h^{-1}x}$, one has
\begin{align}\label{3.22}
 {\mathcal H}_G(x',Z')&=\int_Gc(g^{-1}hh^{-1}x')^2
 {\mathcal H}(g^{-1}x',g^*Z')dg\\
& \geq {1\over 32}c(h^{-1}x)^2{\rm vol}(G_{h^{-1}x})C_{h^{-1}x}|Z'|^2.
\end{align}

By using again the fact that $\overline{U'}$  is compact, one can
cover $U'\cap \mu^{-1}(0)$ by finite open subsets of $M$ verifying
(\ref{3.22}), from which one sees that there exist $C>0$ and
$\varepsilon'>0$ such that for any $x\in U'\cap\mu^{-1}(0)$ and
$Z\in N_x$ with $|Z|\leq \varepsilon'$, one has
\begin{align}\label{3.23}
 {\mathcal H}_G(x,Z) \geq   C|Z|^2.
\end{align}

From (\ref{3.23}) and the $G$-invariance of ${\mathcal H}_G$,
Lemma \ref{t3.5} follows. \ \ Q.E.D.

$\ $

From Lemma \ref{t3.5}, one deduces the following key property.

\begin{lemma}\label{t3.6} There exist $\varepsilon_0>0$ and $C>0$ such
that for any $x\in \mu^{-1}(0)$ and  $Z\in N_x$ with $|Z|\leq
\varepsilon_0$, one has
\begin{align}\label{3.24}
 \left|X^{\mathcal H}_G(x,Z)\right| \geq   C|Z|.
\end{align}
\end{lemma}

{\it Proof}. Recall that for any $z\in M$,
\begin{align}\label{3.25}
 X^{\mathcal H}_G(z)=\int_Gc(g^{-1}z)^2X_g^{\mathcal H}dg.
\end{align}

From (\ref{3.25}) and \cite[(1.14)]{TZ98}, one verifies that,
\begin{align}\label{3.26}
(d{\mathcal H}_G)^*(z) &=\left(d\int_Gc(g^{-1}z)^2{\mathcal
H}(g^{-1}z)dg\right)^*\\
&=JX^{\mathcal H}_G(z) +2\int_G {\mathcal
H}(g^{-1}z)c(g^{-1}z)\left((g^*(dc))(z)\right)^*dg.
\end{align}
Thus, one has
\begin{align}\label{3.27}
 X^{\mathcal H}_G(z)=-J(d{\mathcal H}_G)^*(z)+2J\int_G {\mathcal
H}(g^{-1}z)c(g^{-1}z)\left((g^*(dc))(z)\right)^*dg.
\end{align}

From (\ref{3.23}), one sees that there exists $C'>0$ such that
when $z=(x,Z)$ is close enough to $\mu^{-1}(0)$, one has
\begin{align}\label{3.28}
\left|(d{\mathcal H}_G)^*(z)\right|\geq C'|Z|.
\end{align}

On the other hand, since ${\mathcal H}(\mu^{-1}(0))=0$ and
$(d{\mathcal H})|_{\mu^{-1}(0)}=0$, one verifies easily that
\begin{align}\label{3.29}
{\partial \over \partial Z}\left.\left(\int_G {\mathcal
H}(g^{-1}z)c(g^{-1}z)\left((g^*(dc))(z)\right)^*dg\right)\right|_{Z=0}=0.
\end{align}

 From (\ref{3.27})-(\ref{3.29}), the $G$-invariance property, as
 well as the assumption that $M/G$ is compact, one gets
 (\ref{3.24}).\ \ \ Q.E.D.

 $\ $

Formula (\ref{3.24}) is a direct analogue of \cite[(3.17)]{TZ98}
and \cite[Propo. 8.14]{BL91}. By this and by Proposition
\ref{t3.1}, one sees that one can proceed in exactly the same way
as in \cite[Sections 8 and 9]{BL91} and \cite[Section 3]{TZ98} to
prove that for $p>0$ verifying Proposition \ref{t3.1}, when $T>0$
is large enough, one has
\begin{align}\label{3.30}
{\rm ind}\left(P_fD_{+,T}^{L^p}\right)={\rm
ind}\left(D_+^{L_G^p}\right).
\end{align}
Indeed, all one need is to modify suitably according to the
appearance of the cut-off function $f$. And it is easy to see that
this only causes a modification of adding a compact operator to
the Fredholm operators involved and thus does not alter the
indices  in due course.

From (\ref{3.30}), the obvious invariance  of the independence of
${\rm ind}(P_fD_{+,T}^{L^p})$ with respect to $T$ (which follows
from the Fredholm properties) and Definition \ref{t2.5}, one
completes the proof of Theorem \ref{t1.3}.\ \ \ Q.E.D.

\section{Vanishing properties of cokernels for large
$p$}\label{s4}

We take $E=\Lambda^{0,*}((T^*M)\otimes L^p)$ as in \S2 and \S 3.

Following Remark \ref{t2.4}, let $$\widetilde{D}^{L^p}_+:{\bf
H}^1_f(M,\Lambda^{0,\rm even}((T^*M)\otimes L^p))^G\rightarrow
{\bf H}^0_f(\Lambda^{0,\rm odd}((T^*M)\otimes L^p))^G$$ be the
operator defined by
\begin{align}\label{4.1}
\widetilde{D}^{L^p}_+(fs)=fD^{L^p}_+s.
\end{align}

From (\ref{4.1}) and (\ref{2.9}), one verifies that for any $s\in
\Gamma(\Lambda^{0,\rm even}((T^*M)\otimes L^p))^G$,
\begin{align}\label{4.2}
\widetilde{D}^{L^p}_+(fs)=P_fD^{L^p}_+(fs)-P_f(c(df)s).
\end{align}

By (\ref{4.2}) and Proposition \ref{t2.1}, one finds
\begin{align}\label{4.3}
{\rm ind}\left(\widetilde{D}^{L^p}_+\right)={\rm
ind}\left(P_fD_+^{L^p}\right).
\end{align}

From (\ref{2.14}) and (\ref{4.2}), one finds that for any $s'\in
\Gamma(\Lambda^{0,\rm odd}((T^*M)\otimes L^p))^G$,
\begin{multline}\label{4.4}
\left(\widetilde{D}^{L^p}_+\right)^*(fs')=P_fD^{L^p}_-(fs')+P_f(c(df)s')
=fD_-^{L^p}s'+2P_f(c(df)s')\\
=D_-^{L^p}(fs')-c(df)s'+2P_f(c(df)s').
\end{multline}

From (\ref{2.6}) and (\ref{4.4}), one finds that there exists
$C_1>0$ such that
\begin{align}\label{4.5}
\left\|\left(\widetilde{D}^{L^p}_+\right)^*(fs')\right\|_0\geq
 \left\|D_-^{L^p}(fs')\right\|_0-C_1\left\|fs'\right\|_0 .
\end{align}

On the other hand, since $U'$ has compact closure and $f$ has
compact support in $U'$, by proceeding in exactly the same way as
in \cite[Section 2]{MM}, one sees that there exist $C_2,\, C_3>0$
such that for any  $s'\in \Gamma(\Lambda^{0,\rm odd}((T^*M)\otimes
L^p))^G$,
\begin{align}\label{4.6}
  \left\|D_-^{L^p}(fs')\right\|_0^2\geq \left(C_2p-C_3\right)\left\|fs'\right\|_0^2
\end{align}

From (\ref{4.5}) and (\ref{4.6}), one sees that when $p\geq
{2C_3\over C_2}$, one has
\begin{align}\label{4.7}
\left\|\left(\widetilde{D}^{L^p}_+\right)^*(fs')\right\|_0\geq
 \left(\sqrt{C_2p\over 2} -C_1\right)\left\|fs'\right\|_0 .
\end{align}
From (\ref{4.7}), one sees that when $p\geq {\rm max}\{
{8C_1^2\over C_2}, {2C_3\over C_2}\}$, one has
\begin{align}\label{4.8}
\left\|\left(\widetilde{D}^{L^p}_+\right)^*(fs')\right\|_0\geq
 {1\over 2}\sqrt{C_2p\over 2} \left\|fs'\right\|_0 ,
\end{align}
from which one gets
\begin{align}\label{4.9}
 \ker\left(\widetilde{D}^{L^p}_+\right)^*=0 .
\end{align}

From (\ref{4.1}) and (\ref{4.9}), one  deduces the following result.

\begin{thm}\label{t4.1} There exists $p_0\geq 0$ such that for any
integer $p\geq p_0$,
\begin{align}\label{4.10}
{\rm ind}\left(\widetilde{D}^{L^p}_+\right)=\dim\left(\Ker
D^{L^p}_+\right)^G.
\end{align}
\end{thm}

{\it Proof of Theorem \ref{t1.6}}. By the vanishing theorem due to
Borthwick-Uribe \cite{BU}, Braverman \cite{B} and Ma-Marinescu
\cite{MM}, one sees that when $p>0$ is large enough,
\begin{align}\label{4.11}
{\rm ind}\left( {D}^{L^p_G}_+\right)=\dim\left(\Ker
D^{L^p_G}_+\right).
\end{align}

From Theorem \ref{t1.3}, Definition \ref{t2.5}, (\ref{4.3}),
(\ref{4.10}) and (\ref{4.11}), one sees that the first equality in
(\ref{1.3}) holds when $p>0$ is large enough.

By using (\ref{4.6}), one can proceed as in the proof of
(\ref{4.8}) to see that when $p>0$ large enough, the second
equality in (\ref{1.3}) also holds.

The proof of Theorem \ref{t1.6} is completed. \ \ \ Q.E.D.

\begin{rem}\label{t4.2} Formula (\ref{4.9}) and the second equality in (\ref{1.3})
might be regarded as  extensions of the vanishing theorem of the
half kernel of Spin$^c$-Dirac operators due to Borthwick-Uribe
\cite{BU}, Braverman \cite{B} and Ma-Marinescu \cite{MM} to the
noncompact case.
\end{rem}

\section{Examples and applications}

The main source of examples can be found in the papers of Hochs-Landsman \cite{HL}
and Hochs \cite{H, Hochs-thesis}. In some of their examples, zero is not in the image of the moment 
map.  In this case, by Theorem \ref{t1.6}, for $p$ sufficiently large, we deduce that 
the $G$-invariant kernel of $D^{L^p}_+$ and of $D^{L^p}_-$ both vanish. That is, 
for $p$ sufficiently large,  $D^{L^p}$ is invertible on the $G$-invariant sections in this case.

Another collection of examples, consists of a finitely generated discrete group $G$
acting properly on a locally compact symplectic manifold $(M, \omega)$ such that $M/G$ is compact.
In this case, the moment map is trivial, so that the symplectic reduction of $M$ is just $M/G$, 
which is generally only an orbifold. 
In the special case when $G$ acts freely and properly discontinuously on $M$, so that $M/G$ is a manifold, then Theorem  \ref{t1.3} is well known and for instance it
can be deduced from a result of Pierrot \cite{pierrot}.

\section*{Acknowledgments} Part of the work was done
while the second author   was visiting University of Adelaide in
January of 2008. The work of the first author was partially
supported by the Australian Research Council. The work of the
second author was partially supported by the National Natural
Science Foundation of China.

\newcommand{\diag}{{\rm diag}}
 \newcommand{\dist}{{\rm dist}}
\newcommand{\kaaa}{{\mathfrak k}}
\newcommand{\paaa}{{\mathfrak p}}
\newcommand{\vp}{{\varphi}}
\newcommand{\taaa}{{\mathfrak t}}
\newcommand{\haaa}{{\mathfrak h}}
\newcommand{\R}{{\mathbb R}}
\newcommand{\Hh}{{\bf H}}
\newcommand{\Q}{{\mathbb{Q}}}
\newcommand{\str}{{\rm str}}
\newcommand{\triv}{{\rm triv}}
\newcommand{\Z}{{\mathbb{Z}}}
\newcommand{\bD}{{\bf D}}
\newcommand{\bF}{{\bf F}}
\newcommand{\tX}{{\tt X}}
\newcommand{\Cliff}{{\rm Cliff}}
\newcommand{\tY}{{\tt Y}}
\newcommand{\tZ}{{\tt Z}}
\newcommand{\tV}{{\tt V}}
\newcommand{\tR}{{\tt R}}
\newcommand{\Fam}{{\rm Fam}}
\newcommand{\Cusp}{{\rm Cusp}}
\newcommand{\bK}{{\bf K}}
\newcommand{\K}{{\mathbb{K}}}
\newcommand{\tH}{{\tt H}}
\newcommand{\bS}{{\bf S}}
\newcommand{\bB}{{\bf B}}
\newcommand{\tW}{{\tt W}}
\newcommand{\tF}{{\tt F}}
\newcommand{\bA}{{\bf A}}
\newcommand{\bL}{{\bf L}}
 \newcommand{\bom}{{\bf \Omega}}
\newcommand{\ve}{{\varepsilon}}
\newcommand{\C}{{\mathbb{C}}}
\newcommand{\gen}{{\rm gen}}

\newcommand{\bP}{{\bf P}}
\newcommand{\Naaa}{{\bf N}}
\newcommand{\image}{{\rm image}}
\newcommand{\gaaa}{{\mathfrak g}}
\newcommand{\zaaa}{{\mathfrak z}}
\newcommand{\saaa}{{\mathfrak s}}
\newcommand{\laaa}{{\mathfrak l}}
\newcommand{\stimes}{{\times\hspace{-1mm}\bf |}}
\newcommand{\ausg}{{\rm end}}
\newcommand{\bff}{{\bf f}}
\newcommand{\maaa}{{\mathfrak m}}
\newcommand{\aaaa}{{\mathfrak a}}
\newcommand{\naaa}{{\mathfrak n}}
\newcommand{\brr}{{\bf r}}
\newcommand{\res}{{\rm res}}
\newcommand{\Aut}{{\rm Aut}}
\newcommand{\Pol}{{\rm Pol}}
\newcommand{\Tr}{{\rm Tr}}
\newcommand{\cT}{{\mathcal T}}
\newcommand{\dom}{{\rm dom}}
\newcommand{\db}{{\bar{\partial}}}
\newcommand{\g}{{\gaaa}}
\newcommand{\cZ}{{\mathcal Z}}
\newcommand{\cH}{{\mathcal H}}
\newcommand{\cM}{{\mathcal M}}
\newcommand{\interi}{{\rm int}}
\newcommand{\singsupp}{{\rm singsupp}}
\newcommand{\cE}{{\mathcal E}}
\newcommand{\ccR}{{\mathcal R}}
\newcommand{\cV}{{\mathcal V}}
\newcommand{\cY}{{\mathcal Y}}
\newcommand{\cW}{{\mathcal W}}
\newcommand{\cI}{{\mathcal I}}
\newcommand{\cK}{{\mathcal K}}
\newcommand{\cA}{{\mathcal A}}
\newcommand{\cEp}{{{\mathcal E}^\prime}}
\newcommand{\Ext}{{\mbox{\rm Ext}}}
\newcommand{\rk}{{\rm rank}}
\newcommand{\im}{{\mbox{\rm im}}}
\newcommand{\sign}{{\rm sign}}
\newcommand{\spann}{{\mbox{\rm span}}}
\newcommand{\symm}{{\mbox{\rm symm}}}
\newcommand{\cF}{{\mathcal F}}
\newcommand{\cD}{{\mathcal D}}
\newcommand{\Ree}{{\rm Re }}
\newcommand{\Res}{{\mbox{\rm Res}}}
\newcommand{\Imm}{{\rm Im}}
\newcommand{\inter}{{\rm int}}
\newcommand{\clo}{{\rm clo}}
\newcommand{\Li}{{\rm Li}}
\newcommand{\cN}{{\mathcal N}}
 \newcommand{\conv}{{\rm conv}}
\newcommand{\op}{{\mbox{\rm Op}}}
\newcommand{\cs}{{c_\sigma}}
\newcommand{\ctg}{{\rm ctg}}
\newcommand{\degg}{{\mbox{\rm deg}}}
\newcommand{\Ad}{{\mbox{\rm Ad}}}
\newcommand{\ad}{{\mbox{\rm ad}}}
\newcommand{\codim}{{\rm codim}}
\newcommand{\Gr}{{\mathrm{Gr}}}
\newcommand{\coker}{{\rm coker}}
\newcommand{\id}{{\mbox{\rm id}}}
\newcommand{\ord}{{\rm ord}}
\newcommand{\nat}{{\Bbb  N}}
\newcommand{\sing}{{\mbox{\rm sing}}}
\newcommand{\Ann}{{\mbox{\rm Ann}}}
\newcommand{\aca}{{\aaaa_\C^\ast}}
\newcommand{\acag}{{\aaaa_{\C,good}^\ast}}
\newcommand{\acage}{{\aaaa_{\C,good}^{\ast,extended}}}
\newcommand{\tck}{{\tilde{\ck}}}
\newcommand{\tnk}{{\tilde{\ck}_0}}
\newcommand{\ceep}{{{\mathcal E}(E)^\prime}}
 \newcommand{\ncE}{{{}^\naaa\cE}}
 \newcommand{\Or}{{\rm Or }}
\newcommand{\Diff}{{\mathcal D}iff}
\newcommand{\cB}{{\mathcal B}}
\newcommand{\hc}{{{\mathcal HC}(\gaaa,K)}}
\newcommand{\hcma}{{{\mathcal HC}(\maaa_P\oplus\aaaa_P,K_P)}}
\def\imath{{\rm i}}
\newcommand{\vsl}{{V_{\sigma_\lambda}}}
\newcommand{\czg}{{\cZ(\gaaa)}}
\newcommand{\csl}{{c_{\sigma,\lambda}}}
\def\hB{\hspace*{\fill}$\Box$ \newline\noindent}
\newcommand{\varho}{\varrho}
\newcommand{\Indu}{{\rm Ind}}
\newcommand{\Fin}{{\mbox{\rm Fin}}}
\newcommand{\cS}{{S}}
\newcommand{\orig}{{\mathcal O}}
\def\hB{\hspace*{\fill}$\Box$ \\[0.5cm]\noindent}
 \newcommand{\cG}{{\mathcal G}}
\newcommand{\npci}{{\naaa_P\hspace{-1.5mm}-\hspace{-2mm}\mbox{\rm coinv}}}
\newcommand{\pki}{{(\paaa,K_P)\hspace{-1.5mm}-\hspace{-2mm}\mbox{\rm inv}}}
\newcommand{\mki}{{(\maaa_P\oplus \aaaa_P, K_P)\hspace{-1.5mm}-\hspace{-2mm}\mbox{\rm inv}}}
\newcommand{\Mat}{{\rm Mat}}
\newcommand{\npi}{{\naaa_P\hspace{-1.5mm}-\hspace{-2mm}\mbox{\rm inv}}}
\newcommand{\ngp}{{N_\Gamma(\pi)}}
\newcommand{\gbg}{{\Gamma\backslash G}}
\newcommand{\gkm}{{ Mod(\gaaa,K) }}
\newcommand{\ggkm}{{  (\gaaa,K) }}
\newcommand{\pkm}{{ Mod(\paaa,K_P)}}
\newcommand{\ppkm}{{  (\paaa,K_P)}}
\newcommand{\makm}{{Mod(\maaa_P\oplus\aaaa_P,K_P)}}
\newcommand{\mmakm}{{ (\maaa_P\oplus\aaaa_P,K_P)}}
\newcommand{\cP}{{\mathcal P}}
\newcommand{\gm}{{Mod(G)}}
\newcommand{\gk}{{\Gamma_K}}
\newcommand{\La}{{\mathcal L}}
\newcommand{\cug}{{\cU(\gaaa)}}
\newcommand{\cuk}{{\cU(\kaaa)}}
\newcommand{\dc}{{C^{-\infty}_c(G) }}
\newcommand{\gdk}{{\gaaa/\kaaa}}
\newcommand{\dgkm}{{ D^+(\gaaa,K)-\mbox{\rm mod}}}
\newcommand{\dgm}{{D^+G-\mbox{\rm mod}}}
\newcommand{\vect}{{\C-\mbox{\rm vect}}}
 \newcommand{\cig}{{C^{ \infty}(G)_{K} }}
\newcommand{\gami}{{\Gamma\hspace{-1.5mm}-\hspace{-2mm}\mbox{\rm inv}}}
\newcommand{\cQ}{{\mathcal Q}}
\newcommand{\mmap}{{Mod(M_PA_P)}}
\newcommand{\bbbz}{{\bf Z}}
 \newcommand{\cX}{{\mathcal X}}
\newcommand{\bH}{{\bf H}}
\newcommand{\pr}{{\rm pr}}
\newcommand{\bX}{{\bf X}}
\newcommand{\bY}{{\bf Y}}
\newcommand{\bV}{{\bf V}}

\begin{center}

\section*{Appendix: 
A $KK$-theoretic interpretation of the $G$-invariant index of
Mathai-Zhang} 

by\\ Ulrich Bunke\\
 NWF I - Mathematik, Universit\"at Regensburg, 93040 Regensburg, Germany,\\
email: ulrich.bunke@mathematik.uni-regensburg.de

 \end{center}

\appendix
\section{The class $[D]$}

We consider a locally compact group $G$ which acts properly on a manifold $M$ with compact quotient. On $M$ we consider a $G$-invariant Riemannian metric and a $G$-invariant generalized Dirac operator
$D$ which acts on sections of a $G$-equivariant bundle $F\to M$ equipped with a $G$-equivariant Dirac bundle structure.
The $G$-$C^*$-algebra $C_0(M)$ acts on the $G$-Hilbert space
$\cE:=L^2(M,F)$. The Riemannian manifold $M$ is complete, and therefore by  \cite{MR0369890} the operator $D$ is essentially selfadjoint on
this Hilbert space with domain $C_c(M,F)$.
We consider the $G$-invariant operator
$\cF:=D(D^2+1)^{-1/2}$ defined by function calculus applied to the unique selfadjoint extension  of $D$. The pair
$(L^2(M,F),\cF)$ is a $G$-equivariant Kasparov module (see \cite{MR918241},  \cite{MR1656031}) over $(C_0(M),\C)$ and represents a class $[D]\in KK^G(C_0(M),\C)$.

The group
$$K^G(M):=KK^G(C_0(M),\C)$$
is called the $G$ equivariant $K$-homology group of $M$.

\section{Descent and assembly}\label{uii}

Let $C^*(G)$ denote the maximal group $C^*$-algebra.
In this subsection we provide an explicit description
of the Baum-Connes  assembly map
$$\mu:K^G(M)\to K(C^*(G))\ .$$
Let us first fix some conventions.
If $A$ is a $C^*$-algebra with an action $\rho:G\to \Aut(A)$, then we define
the convolution product on
$C_0(G,A)$ by
$$\phi*\psi(h)=\int_G \phi(g) \rho(g)[\psi(h^{-1}g)]dg\ ,$$
where $dg$ denotes the left-invariant Haar measure on $G$.
The modular character
$\delta:G\to \R^*$ is defined by
$$d(g^{-1})=\delta(g) dg\ .$$
The adjoint $*_*:A\to A$ is the anti-involution given by
$$\phi^{*_*}(g):=\rho(g)[\phi(g^{-1})^*]\ ,
$$ where $a\mapsto a^*$ is the anti-involution of $A$.
The $C^*$-algebra $C^*(G,A)$ is the maximal cross product
of $G$ with $A$
and defined as the closure of the convolution algebra
$C_c(G,A)$ with respect to the norm $$\|\phi\|=\sup_{\kappa}\|\kappa(\phi)\|\ ,$$
where the supremum is taken over all $*$-representations $\kappa$ of $C_c(G,A)$.
The maximal group $C^*$-algebra $C^*(G)=C^*(G,\C)$ is obtained in the special case where
$A:=\C$ has the trivial action of $G$.
The other important example for the present note is
 the $C^*$-algebra $C_0(M)$ with the action
$(g,f)\mapsto g^*f$, $g\in G$, $f\in C_0(M)$.

There is the descent homomorphism
 $$j^G:K^G(M)\cong KK^G(C_0(M),\C)\rightarrow KK(C^*(G,C_0(M)),C^*(G))$$
introduced in \cite[3.11]{MR918241}  (we will give the explicit
description in the proof of Lemma \ref{cascsacacsa}). Following
\cite[Ch.10]{MR1711324} we choose a non-negative cut-off function
$c\in C^\infty_c(M)$ such that $\int_G g^*c^2 dg \equiv 1$. Then we
define the projection $P\in C^*(G,C_0(M))$ by
\begin{equation}\label{dekdelded} P(g)=c g^*c \delta(g)^{1/2}\in
C_c(M)\ .
\end{equation}
Since $G$ acts properly on $M$ we observe that
$P\in C_c(G,C_0(M))$. The relations $P^2=P=P^*$ are straightforward to check,
\begin{align*}
P^2(g,m)& =\int_G c(h^{-1}m)c(m)\delta(h)^{1/2} c(g^{-1}hh^{-1}m) c(h^{-1}m)\delta(h^{-1}g)^{1/2}dh\\
&=c(g^{-1}m)c(m)\delta(g)^{1/2}=P(g,m), \\
P^{*}(g,m)&=c(g^{-1}m)c(gg^{-1}m) \delta(g^{-1})^{1/2}\delta(g)=P(g,m).
\end{align*}
Let $[P]\in K_0(C^*(G,C_0(M))\cong KK(\C, C^*(G,C_0(M)))$ be the class
induced by $P$, which is independent of the choice of $c$, since any two such functions $c_0,c_1$  can be joined by a path $c_t:=\sqrt{tc_0+(1-t) c_1^2}$ which induces a corresponding path of projections.

\begin{defn}
The assembly map $\mu:K^G(M)\to K(C^*(G))\cong KK(\C,C^*(G)) $ is defined
as the composition
 $$KK^G(C_0(M),\C)\stackrel{j^G}{\to}KK(C^*(G,C_0(M),C^*(G))\stackrel{[P]\otimes_{C^*(G,C_0(M))}\dots}{\longrightarrow}   KK(\C,C^*(G))\ .$$
\end{defn}

\section{The index}

The non-reduced group $C^*$-algebra of $G$ has the
universal property that any unitary representation of $G$ extends to a
representation of $C^*(G)$.   In particular, the trivial representation of $G$ on $\C$  has an extension
$1:C^*(G)\rightarrow \C$. On the level of $K$-theory it induces a
homomorphism
$I:K_0(C^*(G))\rightarrow K_0(\C)\cong \Z$.

If we identify
$$K(C^*(G))\cong KK(\C,C^*(G))\ ,$$ then
the homomorphism $I$ can be written as a Kapsarov
product $\dots \otimes_{C^*(G)} [1]$, where $[1]\in KK(C^*(G),\C)$
is represented by the Kasparov module $(\C,0)$.

\begin{defn}\label{aassww}
We define
$\ind:K^G(M)\rightarrow \Z$
to be the composition
$$K^G_0(M)\stackrel{\mu}{\rightarrow} K(C^*(G))
\stackrel{I }{\rightarrow} \Z\ .$$
\end{defn}

\section{A model}
Let $\epsilon:G\to \C^*$ be a character.
By  $L^2_{loc}(M,F)_{\epsilon}^G$ we denote the space of locally
square integrable sections of $F$ which transform under $G$ with character $\epsilon$, i.e. which satisfy
$$ g^*\phi=\epsilon(g)\phi$$
for all $g\in G$.
 Multiplication by $c$ defines a map
$c:L^2_{loc}(M,F)_{\epsilon}^G\to L^2(M,F)$.
This map is actually injective and has a closed range
$H_\epsilon\subseteq L^2(M,F)$. In order to see this we define the continuous maps
$E_\epsilon:L^2(M,F)\to L^2_{loc}(M,F)_{\epsilon}^G$ and
$Q_\epsilon:L^2(M,F)\to H_\epsilon$ by
$$E_\epsilon(\phi) :=\int_G \epsilon^{-1}(g)g^*(c \phi) dg, \qquad Q_\epsilon(\phi):=c E_\epsilon(\phi),
$$
\begin{eqnarray*}
h^*E_\epsilon(\phi)(l)&=&\int_G \epsilon^{-1}(g) (c\phi)(g^{-1}h^{-1}l)dg\\
&=&\int_G \epsilon^{-1}(h^{-1}z)(c\phi)(z^{-1}l)dz\\
&=&\epsilon(h)E_\epsilon(\phi)(l)
\end{eqnarray*}

For $\phi\in L^2_{loc}(M,F)_{\epsilon}^G$ we have
$E_\epsilon(c \phi) =\phi$, since
\begin{eqnarray*}
E_\epsilon(c\phi)(h)&=&\int_G \epsilon^{-1}(g)g^*c^2(h) g^*\phi dg\\&=& \int_G g^*c^2(h) \phi(h) dg\\&=&\phi(h).
\end{eqnarray*}
This implies injectivity of $c$, and furthermore  $Q_\epsilon(c\phi)=c \phi$. Therefore
  $Q_\epsilon$ is a projection onto $H_\epsilon$ which is in fact orthogonal for $\epsilon=\delta^{-1/2}$, the square root of the modular character,
   \begin{eqnarray*}
\langle Q_{\delta^{-1/2}}\phi,\psi\rangle &=& \langle c\int_Gg^*(c \phi) \delta^{1/2}(g)dg,\psi\rangle
=\int_G\langle g^*(c\phi),c \psi\rangle  \delta^{1/2}(g)dg\\
&=&\int_G \langle c\phi,(g^{-1})^*(c \psi)\rangle\delta^{1/2}(g)dg
\stackrel{!}{=}\int_G \langle c\phi,g^*(c \psi)\rangle\delta^{1/2}(g)dg\\
&=& \langle \phi,c\int_Gg^*(c \psi)\delta^{1/2}(g)dg\rangle
=\langle \phi,Q_{\delta^{-1/2}}\psi\rangle\ ,
\end{eqnarray*}
where we use for the marked equality that the measure $\delta^{1/2}(g)dg$ is invariant with respect to inversion $g\mapsto g^{-1}$.

The map $c$ further induces an injective map of Sobolev spaces
$$c:H^1_{loc}(M,F)^G_{\epsilon}\to H^1(M,F)\ ,$$ and we let
$H^1_\epsilon\subseteq H^1(M,F)$ denote the closed image under $c$.

We define the operator
$\tilde D_\epsilon:H^1_\epsilon\to H_\epsilon$ by
\begin{equation}\label{uidwqdwqd}
\tilde D_\epsilon (c f)=Q_\epsilon D  (c f)\ .
\end{equation}

\begin{lemma}
The operator $\tilde D_\epsilon:H^1_\epsilon\to H_\epsilon$ is Fredholm.
\end{lemma}
\proof
To see this we first choose a $G$-invariant parametrix $R$ for $D$ which is a $G$-invariant
pseudo-differential operator of order $-1$ with finite propagation.
This operator induces a $G$-equivariant map
$R:L^2_{loc}(M,F)\to H^1_{loc}(M,F)$, and therefore
$$\tilde R_\epsilon: c R E_{\epsilon|H}:H_\epsilon\to H^1_\epsilon\ $$
We have
$$\tilde D_\epsilon\tilde R_\epsilon=Q_\epsilon Dc RE=Q_\epsilon c E_\epsilon+ Q_\epsilon[D,c]RE+Q_\epsilon c (DR-1)E_\epsilon\ .$$
Using that $Q_\epsilon c E_{\epsilon|H_\epsilon}=1_{H_\epsilon}$, and that $[D,c] =c(dc)$ and $c$ are  compactly supported
operators of order zero, and $R$ and $DR-1$ are of order $-1$ we see that  $\tilde D_\epsilon\tilde R_\epsilon-1_{H_\epsilon}$ is compact. In a similar manner we show that $\tilde R_\epsilon\tilde D_\epsilon-1_{H^1_\epsilon}$ is compact.
\hB
Note that the index of the operator $\tilde D_1$ associated to the trivial character is studied in
the main text by Mathai-Zhang. On the other hand, we will see in Proposition \ref{MAIN} that the index of the operator $\tilde D_{\delta^{-1/2}}$ is equal to $\ind([D])$.
The following result connects both cases.

\begin{lemma}
We have $\ind(\tilde D_1)=\ind(\tilde D_{\delta^{-1/2}})$.
\end{lemma}
\proof The main idea is that $\delta^{1/2}$ can be connected with
the trivial character by the continuous path of characters
$\epsilon_t:=\delta^{-t/2}$, $t\in [0,1]$. Let $I=[0,1]$ and
consider the $C(I)$-Hilbert-modules $C(I,L^2(M,F))$ and
$C(I,H^1(M,F))$. The family of operators $Q_{\epsilon_t}$ defines
projections $Q$ on $C(I,L^2(M,F))$ and $C(I,H^1(M,F))$ with images
$H$ and $H^1$. Furthermore, the family $\tilde D_{\epsilon_t}$
induces an operator $\tilde D: H^1\to H^0$ whose parametrix $\tilde
R$ is given by the family $\tilde R_{\epsilon_t}$. These data give a
Kasparov module $$(H\oplus H^1,\begin{pmatrix} 0&\tilde R\\\tilde
D&0\end{pmatrix})$$ over $C(I)$ which is a homotopy (see
\cite{MR1656031}) between the Kasparov modules
$$(H_1\oplus H^1_1,\begin{pmatrix} 0&\tilde R_1\\\tilde D_1&0\end{pmatrix})\ ,\quad (H_{\delta^{-1/2}}\oplus H^1_{\delta^{-1/2}},\begin{pmatrix} 0&\tilde R_{\delta^{-1/2}}\\\tilde D_{\delta^{-1/2}}&0\end{pmatrix})$$
representing $\ind(\tilde D_1)$ and $\ind(D_{\delta^{-1/2}})$.
\hB
We can now formulate the main assertion of this note:

\begin{prop}\label{MAIN}
We have
$\ind([D])=\ind(\tilde D_{\delta^{-1/2}})$
\end{prop}

\section{Proof of \ref{MAIN}}

We first apply $j^G$ to the Kasparov module $(L^2(M,F),\cF)$ representing
$[D]$. Note that by the universal property of the maximal crossed product the compatible  $G$ and $C_0(M)$-actions on $L^2(M,F)$
extend to an action
of $C^*(G,C_0(M))$.
\begin{lemma}\label{cascsacacsa}
$j^G([D])\otimes_{C^*(G)} [1]\in KK(C^*(G,C_0(M)),\C)$ is represented by the Kasparov module
$(L^2(M,F),\cF)$.
\end{lemma}
\proof According to \cite[3.11]{MR918241}, $j^G([D])$ is represented
by $(C^*(G,L^2(M,F)),\tilde\cF)$, where $C^*(G,L^2(M,F))$ is a right
$C^*(G)$-module with a left action by $C^*(G,C_0(M))$. It is a
closure of the space of compactly supported continuous functions
$f:G\rightarrow L^2(M,F)$. The operator $\tilde \cF$ is given by
$(\tilde\cF f)(g)=(\cF f)(g)$. The $C^*(G)$-valued scalar product is
given by $$\langle f_1,f_2\rangle(g)=\int_D\langle f_1(h),
f_2(hg)\rangle  dg\ $$ Furthermore, the left action of
$C^*(G,C_0(M))$ is given by $$(\phi f)(g)=\int_G \phi(h) (hf)(g) dh\
.$$

Using that $C^*(G,L^2(M,F))\otimes _{C^*(G)}
\C\cong L^2(M,F)$ by
$$f\otimes v\mapsto  \int_G  f(g)  dg v\ .$$
Therefore
$j^G([D])\otimes_{C^*(G)} [1]$ is represented by the Kasparov module
$(L^2(M,F),\cF)$, where  the left action of $C^*(G,C_0(M))$ on $L^2(M,F)$ is given  by
\begin{equation}\label{szswqs}
(\phi f)=\int_{G} \phi(h) hf dh\ .
\end{equation}
\hB
We now compute $[P]\otimes_{C^*(G,C_0(M))}\left(j^G([D])\otimes_{C^*(G)}
[1]\right)$. We represent $[P]$ by the Kasparov module
$(PC^*(G,C_0(M)),0)$. We must understand $PC^*(G,C_0(M))
\otimes_{C_0(M)} L^2(M,F)$. The operator $Q_{\delta^{-1/2}}:=c E_{\delta^{-1/2}}$
 is the orthogonal projection
$L^2(M,F)\to H_{\delta^{-1/2}}$.
 If we combine the formula (\ref{szswqs})  for the action of $C^*(G,C_0(M))$ on $L^2(M,F)$ with the definition (\ref{dekdelded}) of $P$ we see that the projection $P$ acts as $Q_{\delta^{-1/2}}$, and
\begin{eqnarray*}
PC^*(G,C_0(M))
\otimes_{C^*(G,C_0(M))}L^2(M,F)&=& P L^2(M,F)\\=Q_{\delta^{-1/2}}L^2(M,F)&=&
H_{\delta^{-1/2}}
\end{eqnarray*}
Let  $L:H_{\delta^{-1/2}}\to L^2(M,F)$ denote the inclusion.
Since $Q_{\delta^{-1/2}}$ is othogonal, the operator (\ref{uidwqdwqd}) can be written in the form
 $$\tilde D=L^*DL$$ which makes clear that it is selfadjoint.   We form $\tilde \cF:=\tilde  D_{\delta^{-1/2}}(1+\tilde D_{\delta^{-1/2}}^*\tilde D_{\delta^{-1/2}})^{-1/2}$.
The Kasparov module $(H,\tilde \cF)$ over $(\C,\C)$ represents
$\ind(\tilde D_{\delta^{-1/2}})\in \Z$ under $KK^0(\C,\C)\cong \Z$.
The assertion of   Proposition \ref{MAIN} immediately follows from

\begin{lemma}
 $[P]\otimes_{C^*(G,C_0(M))}\left(j^G([D])\otimes_{C^*(G)}
[1]\right)$ is represented by the Kapsarov module
$(H,\tilde  \cF)$.
\end{lemma}
\proof In order to show the claim we employ the characterization of
the Kasparov product in terms of connections (see
\cite[2.10]{MR918241}). In our situation we have only to show that
$\tilde \cF$ is an $\cF$-connection.

For Hilbert-$C^*$-modules $X,Y$ over some $C^*$-algebra $A$ let $L(X,Y)$ and $K(X,Y)$ denote the spaces of bounded and compact adjoinable $A$-linear operators (see \cite{MR1656031} for definitions).
For $\xi\in PC^*(G,C_0(M))$ we define $\theta_\xi\in
L(L^2(M,F) , H )$ by $\theta_\xi(f)=\xi f$.
Since $\cF$ and $\tilde \cF$ are selfadjoint,     we only must show that
$$\theta_\xi\circ   \cF- \tilde  \cF \circ \theta_\xi\in
K(L^2(M,F) , H )\ .$$

We have $\xi  \cF- \tilde  \cF \xi
= [\xi, \cF]+( \cF -\tilde  \cF )  P\xi$.
Since $[\xi, \cF]$ is compact it suffices to show that $( \cF -\tilde
\cF )  P$ is compact.
We consider $\bar D:=(1-P)D(1-P)+L\tilde  D L^*$.
Then by a simple calculation we have $\bar D=D+R$, where $R$ is a   bounded operator.
Let $\bar\cF:=\bar\cD(1+\bar \cD^2)^{-1/2}$. Then
$(\cF -L\tilde  \cF L^*) P=(\cF-\bar \cF)P$.
Let $\tilde c\in C_c^\infty(M)$ be such that $c\tilde c=c$.
Then we have $(\hat \cF-\bar\cF)P=(\cF-\bar \cF)\tilde c P$.
Therefore it suffices to show that $(\cF-\bar \cF)\tilde c$ is
compact. This can be done using the integral representations for $ \cF$
and $\bar \cF$ as in \cite{MR1348799}.



\end{document}